\documentclass[thmsa,final,french,11pt]{article}

\usepackage{eurosym}
\usepackage{amsfonts}
\usepackage{amsmath}
\usepackage{amssymb}
\usepackage{mathrsfs} 
\usepackage{enumerate}
\usepackage{graphicx}

\usepackage{ntheorem}

\theoremstyle{definition}
\newtheorem{Def}{Definition}

\newtheorem{Th}{Theorem}

\newtheorem{Prop}[Th]{Proposition}
\newtheorem{Le}[Th]{Lemma}

\newtheorem{Rk}{Remark}

\newcommand{\WW}{\mathbf{W}}

\newcommand{\RR}{\mathbf{R}}
\newcommand{\ZZ}{\mathbf{Z}}

\newcommand{\conv}[2][n]{\underset{#1\rightarrow #2}{\longrightarrow}}
\newcommand{\eq}[2][n]{\underset{#1\rightarrow #2}{\sim}}
\newcommand{\EEE}[1]{\operatorname{\mathbb{E}}\left[\,#1\,\right]}

\newcommand{\PPP}[1]{\operatorname{\mathbb{P}}\left(\,#1\,\right)}

\newcommand{\ind}[1]{\mathbb{I}_{#1}\,}

\newcommand{\sgn}{\text{sgn}}

\def\WW_#1{\boldsymbol{W}\!_#1}
\numberwithin{equation}{section}

\newenvironment{prooft}[1]{\vskip 2mm\noindent {\bf Proof of #1.}}{\hfill $\square$ \noindent}

\usepackage{authblk}
\usepackage{xcolor}

\usepackage[french]{varioref}

\begin{document}
\author[,1]{\textsc{Nicolas Chenavier} \thanks{\texttt{nicolas.chenavier@univ-littoral.fr}}}
\author[,1]{\textsc{Ahmad Darwiche} \thanks{\texttt{ahmad.darwiche@univ-littoral.fr}}}
\author[,2]{\textsc{Arnaud Rousselle} \thanks{\texttt{arnaud.rousselle@u-bourgogne.fr}.}}
 \affil[1]{Universit\'e du Littoral C\^ote d'Opale, UR 2597, LMPA, Laboratoire de Math\'ematiques Pures et Appliqu\'ees Joseph Liouville,
62100 Calais, France.}
\affil[2]{Institut de Math\'ematiques de Bourgogne, UMR 5584, CNRS, Universit\'e Bourgogne Franche-Comt\'e,
 F-21000 Dijon, France.}

\title{Compound Poisson approximation for simple transient random walks in random sceneries}
\maketitle

\begin{abstract}
Given a simple transient random walk $(S_n)_{n\geq 0}$ in $\ZZ$ and a stationary sequence of real random variables $(\xi(s))_{s\in \ZZ}$, we investigate the extremes of the sequence $(\xi(S_n))_{n\geq 0}$. Under suitable conditions, we make explicit the extremal index and show that the point process of exceedances converges to a compound Poisson point process. We give two examples for which the cluster size distribution can be made explicit.
\end{abstract}



\bigskip

\section{Introduction}

Extreme Value Theory (EVT) deals with rare events and has many applications in various domains such as hydrology \cite{KPN}, finance \cite{EKM} and climatology \cite{YGLN}. It was first introduced in the context of independent and identically distributed (i.i.d.) random variables.  It is straightforward that if $(\xi(s))_{s\in \mathbf{Z}}$ is a sequence of i.i.d. random variables then the following property holds: for any sequence of real numbers $(u_n)_{n\geq 0}$, and for $\tau >0$,
\[n\PPP{\xi(0)>u_n}\conv[n]{\infty}\tau \Longrightarrow \PPP{\max_{0\leq k\leq n}\xi(k)\leq u_n}\conv[n]{\infty}e^{-\tau}.\] 
The above property has been extended for sequences of dependent random variables satisfying two conditions. The first one, referred to as the $D(u_n)$ condition of Leadbetter, is a long range dependence property and the second one, known as the $D'(u_n)$ condition, ensures that, locally, there is no clusters of exceedances (see \cite{L2} for a statement of these conditions).

In 2009, Franke and Saigo \cite{franke_saigo2009bis, franke_saigo_2009} considered the following problem. Let $(X_i)_{i\geq 1}$ be a sequence of centered, integer-valued i.i.d. random variables and let $S_0=0$ and $S_{n}=X_{1}+\dots+X_{n}$, $n \geq 1$. Assume that $(X_i)_{i\geq 1}$ is in the domain of attraction of a stable law, i.e. for each $x\in \RR$, 
\begin{equation*}
\PPP{n^{-\frac{1}{\alpha}}S_{n}\leq x}  \conv[n]{\infty}  F_{\alpha}(x),    
\end{equation*} where $F_{\alpha}$ is the distribution function of a stable law with characteristic function given by
\[\varphi(\theta)=\exp(-|\theta|^{\alpha}(C_{1}+iC_{2}\sgn \theta)), \ \alpha \in (0,2].\]  When $\alpha<1$ (resp. $\alpha>1$), it is known that the random walk $(S_n)_{n\geq 0}$ is transient (resp. recurrent) \cite{KS,LeGRo}. Now, let $(\xi(s))_{s\in \ZZ}$ be a family of $\RR$-valued which are i.i.d. random variables independent of the sequence  $(X_i)_{i\geq 0}$. In the sense of \cite{franke_saigo_2009}, the sequence $(\xi(S_n))_{n\geq 0}$ is called a \textit{random walk in a random scenery}. Franke and Saigo derive limit theorems for the maximum of the first $n$ terms of $(\xi(S_n))_{n\geq 0}$ as $n$ goes to infinity. An adaptation of Theorem 1 in \cite{franke_saigo_2009} shows that in the transient case, i.e. $\alpha<1$, the following property holds: if $n\PPP{\xi(0)>u_n}\conv[n]{\infty}\tau$ for some sequence $(u_n)_{n\geq 0}$ and for some $\tau>0$, then \begin{equation}\label{eq:frankesaigo}  \PPP{\max_{0\leq k\leq n}\xi(S_k)\leq u_n} \conv[n]{\infty}e^{-q\tau },
\end{equation} where 
\begin{equation} 
\label{eq:defq}
q=\PPP{S_i\neq 0, \forall i\geq 1}.
\end{equation} Notice that $q>0$ because the random walk $(S_n)_{n\geq 0}$ is transient.  The term $q$ can also be expressed as (see e.g. \cite{LeGRo})
\begin{equation}\label{eq:range} q=\lim_{n\rightarrow \infty}\frac{R_n}{n} \quad \text{a.s.},
\end{equation} where $R_n=\#\{S_0,\ldots, S_n\}$ is the \textit{range} of the random walk.  The result \eqref{eq:frankesaigo} was recently extended  to a random scenery which is not necessarily based on i.i.d. random variables but on a sequence satisfying a slight modification of the $D(u_n)$ condition \cite{Chen_Dar}. One of the difficulties is that the sequence $(\xi(S_n))_{n\geq 0}$ does not satisfy the $D'(u_n)$ condition and clusters of exceedances can appear. 

In this paper, we establish more precise results on the long range dependence and on the so-called point process of exceedances. To do it, we assume some additional properties on the random walk and on the random scenery. First, we let $S_0=0$ and $S_n=X_1+\cdots +X_n$, $n\geq 1$, where  $(X_i)_{i\geq 1}$ is a family of independent random variables with distribution $\PPP{X_i=1}=p$ and $\PPP{X_i=-1}=1-p$, $p\in ]0,1[\setminus \{\tfrac{1}{2}\}$. Notice that, as opposed to  \cite{Chen_Dar} and \cite{franke_saigo_2009}, the random variable $X_i$ is not centered. In what follows, the term $q$ appearing in \eqref{eq:defq} also satisfies \eqref{eq:range}. It can easily be proved that $q=|2p-1|$. Secondly, we assume that the random scenery $(\xi(s))_{s\in \ZZ}$ satisfies a long range dependence property, which is referred to as the $\Delta(u_n)$ condition. To state it, we give some notation.  For each $n,m_1,m_2$ with $0\leq m_1 \leq m_2\leq n$, define $\mathcal{B}_{m_1}^{m_2}(u_n)$ as the $\sigma$-field generated by events $\{\xi(s)\leq u_n\}$, $m_1\leq s\leq m_2$, where $(u_n)_{n\geq 0}$ is some sequence of positive numbers. Also for each $n$ and $1\leq \ell\leq n-1$, write 
 \[\alpha_{n,\ell} =  \max_{1\leq k\leq n-\ell} \max_{A\in \mathcal{B}_0^k(u_n), B\in \mathcal{B}_{k+\ell}^n(u_n) }\{ |\PPP{A\cap B} - \PPP{A}\PPP{B}|\}.\] 
 We are now prepared to state the $\Delta(u_n)$ condition.
\begin{Def}
 We say that the stationary sequence $(\xi(s))_{s\in \ZZ}$ satisfies the $\Delta(u_n)$ condition if there exists some sequence $(\ell_n)_{n\geq 0}$ such that $\alpha_{n,\ell_n}\conv[n]{\infty}0$ and  $\ell_n=o(n)$. 
\end{Def} 
The above condition is slightly more restrictive than the $D(u_n)$ condition and was introduced by Hsing, Hüsler and Leadbetter \cite{HHL} in the context of stationary sequence of random variables indexed by the set of \textit{positive} integers. 

\paragraph{The extremal index}
Assume from now on that, for any $\tau>0$, there exists a threshold $u_n=u_n^{(\tau)}$ such that
 \begin{equation}
 \label{eq:threshold}
 n\PPP{\xi(0)>u_n}\conv[n]{\infty}\tau.
 \end{equation}
The existence of the threshold $u_n=u_n^{(\tau)}$ is ensured when $\lim_{x\rightarrow x_f} \frac{\overline{F}(x)}{\overline{F}(x-)}=1$, where $x_f = \sup\{u: F(u)<1\}$, $F(u)=\PPP{\xi(0)\leq u}$ and $\overline{F}=1-F$ (see Theorem 1.1.13 in \cite{LLR}). 

To state our first main result, we recall the $D^{(k)}(u_n)$ condition, introduced by Chernick, Hsing and McCormick \cite{CHM}.

\begin{Def}
\label{def:dkun}
Let $k\geq 1$. Assume that $(\xi(s))_{s\in \ZZ}$ satisfies the $\Delta(u_n)$ condition. We say that $(\xi(s))_{s\in \ZZ}$ satisfies the $D^{(k)}(u_n)$ condition if there exist sequences of integers $(s_n)_{n\geq 0}$ and $(\ell_n)_{n\geq 0}$ such that $s_n\conv[n]{\infty}\infty$, $s_n \alpha_{n, \ell_n}\conv[n]{\infty}0$, $s_n\ell_n/n\conv[n]{\infty}0$, and
\[\lim_{n\rightarrow \infty}  n\PPP{\xi(0)>u_n\geq M_{1,k-1}, M_{k,r_n}>u_n} = 0, \]
where $M_{i,j}=-\infty$ for $i>j$, $M_{i,j}=\max_{0\leq i\leq t\leq j}\xi(t)$ for $i\leq j$, and $r_n=\lfloor n/s_n\rfloor$.
\end{Def}
As noticed in \cite{CHM}, the $D'(u_n)$ condition is slightly more restrictive than the $D^{(1)}(u_n)$ condition. Observe that the $D^{(k)}(u_n)$ condition is satisfied when the sequence $(\xi(s))_{s\in \ZZ}$ is $k$-dependent. Recall also that the (stationary) sequence $(\xi(s))_{s\in \ZZ}$ has an \textit{extremal index} $\sigma \in [0,1]$ if, in conjunction to \eqref{eq:threshold}, we have
 \[\PPP{\max_{0\leq s\leq n} \xi(s) \leq u_n} \conv[n]{\infty}e^{-\sigma \tau},\]
  for any $\tau >0$.  The extremal index can be interpreted as the reciprocal of the mean size of a cluster of exceedances. According to Corollary 1.3. in \cite{CHM}, under the assumptions that the $\Delta(u_n)$ and $D^{(k)}(u_n)$ conditions hold for $u_n=u_n^{(\tau)}$ for any $\tau >0$, the extremal index exists and is equal to $\sigma$ if and only if $\PPP{M_{1,k}\leq u_n|\xi(0)>u_n}\conv[n]{\infty}\sigma$ for any $\tau >0$.  In particular, when the $D'(u_n)$ condition is satisfied, the extremal index exists and is equal to $\sigma=1$ (Theorem 1.2. in \cite{L2}). 
 
The following proposition ensures that, under suitable conditions, the extremal index of the sequence $(\xi(S_n))_{n\geq 0}$ exists and can be made explicit.
 
 \begin{Prop}
\label{prop:existence_extremal_index}
Let  $(\xi(s))_{s\in \ZZ}$ be a stationary sequence satisfying the $D^{(k)}(u_n)$ condition with $k\geq 1$, for $u_n=u_n^{(\tau)}$ for any $\tau>0$ and such that $(\xi(S_n))_{n\geq 0}$ satisfies the $\Delta (u_n)$ condition. Assume that the extremal index of $(\xi(s))_{s\in \ZZ}$ exists and is equal to $\sigma \in (0,1]$, i.e.
\[n\PPP{\xi(0)>u_n} \conv[n]{\infty}\tau \quad \text{and} \quad \PPP{\max_{0\leq s\leq n} \{\xi(s)\} \leq u_n}\conv[n]{\infty}e^{-\sigma\tau}. \] Then the sequence $(\xi(S_n))_{n\geq 0}$ admits an extremal index which is equal to $\theta=\sigma q$, where $q$ is as in \eqref{eq:defq}.
\end{Prop}
In other words, Proposition \ref{prop:existence_extremal_index} claims that  $\PPP{\max_{0\leq k\leq n}\xi(S_k)\leq u_n} \conv[n]{\infty}e^{-\sigma q \tau }$. 

\paragraph{The $D^{(k)}(u_n)$ condition}
In \cite{franke_saigo_2009}, Franke and Saigo proved that, when the $\xi(s)$'s are i.i.d., the sequence $(\xi(S_n))_{n\geq 0}$ does not satisfy the $D'(u_n)$ condition. Since the $D^{(1)}(u_n)$ condition is slightly less restrictive than the $D'(u_n)$ condition and since Equation \eqref{eq:sufficientconditionnotdk} (as stated below) is satisfied when the $\xi(s)$'s are i.i.d., the following result can be compared to Proposition 3 in \cite{franke_saigo_2009}.

\begin{Prop}
\label{prop:dkun}
Assume that $\left(\xi(s)\right)_{s\in\ZZ}$ is a stationary sequence satisfying
\begin{equation}
\label{eq:sufficientconditionnotdk}
 \limsup_{n\to\infty}n\mathbb{P}\Big({\xi(0)>u_n\mbox{ and there exists }s\in\{-k+1,\dots,k-1\}\setminus\{0\}\mbox{ s.t. }\xi(s)>u_n}\Big)<\tau.
 \end{equation}
 Then $(\xi(S_n))_{n\geq 0}$ does not satisfy the $D^{(k)}(u_n)$ condition for any $k\geq 1$. 
\end{Prop}

\paragraph{The point process of exceedances}
Before stating our next result, we give some notation. First, for any $\tau>0$ and $n\geq 1$, we let
 \begin{equation}
 \label{eq:defppexceedances}
 \Phi_n=\Phi_n^{(\tau)} = \left\{ \frac{i}{n}: \xi(S_i)>u_n, i\leq n  \right\}.
 \end{equation}
 The point process $\Phi_n\subset [0,1)$ is referred to as the \textit{point process of exceedances}.   Now, let $(k_n)_{n\geq 0}$ be a sequence of positive integers with $k_n\conv[n]{\infty}\infty$ and $k_n=o(n)$. For any $j\geq 1$, let
 \[p_n(j) = p_n^{(\tau)}(j)= \PPP{\#\Phi_{k_n} = j | \xi(0)>u_n},\]
 where $\Phi_{k_n}$ is defined in the same spirit as $\Phi_n$ by considering this time integers $i\leq k_n$. Recall that a \textit{compound Poisson point process} in $[0,1]$ of intensity $\lambda>0$ and cluster size distribution $\pi$  is a point process of the form $\{(x_j, n_j), j\geq 1\}$, where $\{x_j, j\geq 1\}$ is a stationary Poisson point process in $[0,1]$ of intensity $\lambda$ and where $(n_j)_{j\geq 1}$ is a family of i.i.d. random variables with distribution $\pi$ which is independent of the $x_j$'s. 
 
The following proposition states that, under suitable conditions, the point process of exceedances converges to a compound Poisson point process. 
 
\begin{Prop}
\label{prop:ppexceedances}

Let  $(\xi(s))_{s\in \ZZ}$ be a stationary sequence satisfying the $D^{(k)}(u_n)$ condition with $k\geq 1$, for $u_n=u_n^{(\tau)}$ for any $\tau>0$ and such that $(\xi(S_n))_{n\geq 0}$ satisfies the $\Delta (u_n)$ condition. Assume that $p_n^{(\tau_0)}(j)$ converges to some number $p(j)$ for any $j\geq 1$ and for some $\tau_0$. Then the point process $\Phi_n^{(\tau)}$ converges in distribution to a compound Poisson point process with intensity $\theta \tau$, where $\theta=\sigma q$, and cluster size distribution $\pi(j)=\frac{1}{\theta}(p(j)-p(j+1))$, for any $\tau>0$.
\end{Prop}
The above proposition is classical in EVT and is a simple application of Theorem 2.5 in \cite{Perfekt94} and Theorem 5.1 in \cite{HHL}. In a different context, asymptotic results on point processes associated with extremes in random sceneries are also established in \cite{Dar2}.

Our paper is organized as follows. In Section \ref{sec:proofs}, we prove Propositions \ref{prop:existence_extremal_index}-\ref{prop:ppexceedances}. In Section \ref{sec:examples}, we give some examples illustrating Proposition \ref{prop:ppexceedances}. In particular, we make explicit the cluster size distribution of the limiting point process of exceedances. In section \ref{sec:discussion}, we shortly discuss possible ways to ensure that $(\xi(S_n))_{n\geq 0}$ satisfies the $\Delta (u_n)$ condition to apply Propositions \ref{prop:existence_extremal_index} and \ref{prop:ppexceedances}. Finally, the proof of a (technical) side-lemma is given in Appendix \ref{app:proofconcentrationinequality}.

\section{Proofs of the main results} 
\label{sec:proofs}

\subsection{Proof of Proposition \ref{prop:existence_extremal_index} }  
Let $R_n=\#\{S_0,\ldots, S_n\}$ be the range associated with the random walk $(S_n)_{n\geq 0}$. Then
\[
\PPP{\max_{0\leq i\leq n}\xi(S_i)\leq u_n} = \EEE{\PPP{\max_{1\leq s\leq R_n}\{\xi(s)\}\leq u_n |R_n  }}.
\]
Moreover,
\begin{align*}
\left| \PPP{\max_{1\leq s\leq R_n}\{\xi(s)\}\leq u_n | R_n}  - \PPP{\max_{1\leq s\leq \lfloor qn\rfloor}\{\xi(s)\}\leq u_n}  \right| & \leq 2\PPP{\exists s \in (R_n,\lfloor qn\rfloor): \xi(s)>u_n}\\
& \leq 2 |R_n-\lfloor qn\rfloor|\PPP{\xi(0)>u_n},
\end{align*}
where $(R_n,\lfloor qn\rfloor)$ denotes the interval with (non-necessarily ordinated) extremities $R_n$ and $\lfloor qn\rfloor$. According to \eqref{eq:range} and \eqref{eq:threshold},  we deduce that \[\PPP{\max_{1\leq s\leq R_n}\{\xi(s)\}\leq u_n | R_n}  - \PPP{\max_{1\leq s\leq \lfloor qn\rfloor}\{\xi(s)\}\leq u_n} \conv[n]{\infty}0.\] Therefore, to prove that $(\xi(S_n))_{n\geq 0}$ has an extremal index which is equal to $\theta=\sigma q$, it is sufficient to prove that 
\begin{equation}
\label{eq:extremalindexxiqn}
 \PPP{\max_{1\leq s\leq \lfloor qn\rfloor}\{\xi(s)\}\leq u_n} \conv[n]{\infty} e^{-\sigma q\tau}.
 \end{equation}

We prove below \eqref{eq:extremalindexxiqn}.

When $k=1$, the identity follows from Corollary 1.3 in \cite{CHM} which also shows that the extremal index is $\sigma=1$. Assume from now on that $k\geq 2$. Because $(\xi(s))_{s\in \ZZ}$ satisfies the $D^{(k)}(u_n)$ condition, it follows from Corollary 1.3 in \cite{CHM} that 
\[\PPP{\max_{1\leq s\leq k-1}\{\xi(s)\}\leq u_n| \xi(0)>u_n }\conv[n]{\infty}\sigma.\] In particular
\[\PPP{\max_{1\leq s\leq k-1}\{\xi(s)\}\leq u_n^{(q)}| \xi(0)>u_n^{(q)} }\conv[n]{\infty}\sigma,\]
where $u_n^{(q)}=u_{\lfloor n/q \rfloor}$. Observe that $(\xi(s))_{s\in \ZZ}$ also satisfies the $\Delta(u^{(q)}_n)$ and $D^{(k)}(u^{(q)}_n)$ condition.  Because $n\PPP{\xi(0)>u_n^{(q)}}\conv[n]{\infty}q\tau $, it follows again from Corollary 1.3 in \cite{CHM} that 
\[ \PPP{\max_{1\leq s\leq n}\{\xi(s)\}\leq u_n^{(q)}} \conv[n]{\infty} e^{-\sigma q\tau}.\]
Taking $n=\lfloor qn'\rfloor$, we deduce that 
\[ \PPP{\max_{1\leq s\leq \lfloor qn'\rfloor }\{\xi(s)\}\leq u^{(q)}_{\lfloor qn'\rfloor}} \conv[n']{\infty} e^{-\sigma q\tau}.\] 
Because  $ \frac{n'}{\lfloor \lfloor qn'\rfloor/q  \rfloor} \conv[n']{\infty}1$, the latter expression gives 
\[ \PPP{\max_{1\leq s\leq \lfloor qn'\rfloor}\{\xi(s)\}\leq u_{n'}} \conv[n']{\infty} e^{-\sigma q\tau},\] 
which proves \eqref{eq:extremalindexxiqn}.

\subsection{Proof of Proposition \ref{prop:dkun} }
It is sufficient to prove that, for any sequence $(r_n)_{n\geq 0}$ with $r_n\conv[n]{\infty}\infty$, we have
\begin{equation}
\label{eq:aimdkun}
\liminf_{n\rightarrow\infty} n\PPP{\xi(S_{{{}}0})>u_n\geq \tilde{M}_{{{}1},{{}k-1}}, \tilde{M}_{{{}k},r_n}>u_n} >0,
\end{equation}
where, similarly to Definition \ref{def:dkun}, we let  $\tilde{M}_{i,j}=-\infty$ for $i>j$, $\tilde{M}_{i,j}=\max_{0\leq i\leq t\leq j}\xi(S_t)$ for $i\leq j$.

Assume that $k$ is {{} even}. We have
\begin{multline}
\label{eq:maintermdkun}
n\PPP{\xi(S_{{}0})>u_n\geq \tilde{M}_{{{}1},{{}k-1}}, \tilde{M}_{{{}k},r_n}>u_n}\\
\begin{split}
& \geq  n\PPP{\xi({{}0})>u_n\geq \tilde{M}_{{{}1},{{}k-1}}, S_{{{}k}}={{}0}}\\
& =  n\PPP{\xi({{}0})>u_n, S_{{{}k}}={{}0}} - n\PPP{\xi({{}0})>u_n, \tilde{M}_{{{}1},{{}k-1}}>u_n, S_{{{}k}}={{}0}}.
\end{split}
\end{multline}
First, because $(\xi(s))_{s\in \ZZ}$ and $(S_n)_{n\geq 0}$ are independent, we obtain from \eqref{eq:threshold} that
\begin{equation}
\label{eq:maintermdkunpart1}
n\PPP{\xi({{}0})>u_n, S_{{{}k}}={{}0}} \eq[n]{\infty} \tau \PPP{S_{{{}k}}={}0},
\end{equation}
where $\PPP{S_{{{}k}}={{}0}}\neq 0$ since $k$ is {{}even}. Secondly, we have
\begin{multline}
\label{eq:prooftermdkun}
n\PPP{\xi({{}0})>u_n, \tilde{M}_{{{}1},{{}k-1}}>u_n, S_{{{}k}}={{}0}}\\
 = n\PPP{\left\{\xi({{}0})>u_n,  S_{{{}k}}={{}0}\right\} \cap \bigcup_{{{}1}\leq i\leq {{}k-1}}\{S_i = {{}0}\} }\\
+ n\PPP{\left\{\xi({{}0})>u_n, \tilde{M}_{{{}1},{{}k-1}}>u_n, S_{{{}k}}={{}0}\right\} \cap \bigcap_{{{}1}\leq i\leq {{}k-1}}\{S_i\neq {{}0}\} }. 
\end{multline}
The first term of the right hand-side of \eqref{eq:prooftermdkun} is equal to 
\[n \PPP{\xi({{}0})>u_n} \PPP{\left\{ S_{{{}k}}={{}0}\right\} \cap \bigcup_{{{}1}\leq i\leq {{}k-1}}\{S_i = {{}0}\} } \eq[n]{\infty} \tau \PPP{\left\{ S_{{{}k}}={{}0}\right\} \cap \bigcup_{{{}1}\leq i\leq {{}k-1}}\{S_i = {{}0}\} }. \]

 To deal with the second term of the right hand-side of \eqref{eq:prooftermdkun}, observe that
 
 \begin{multline*}
 n \PPP{\left\{\xi({{}0})>u_n, \tilde{M}_{{{}1},{{}k-1}}>u_n, S_{{{}k}}={{}0}\right\} \cap \bigcap_{{{}1}\leq i\leq {{}k-1}}\{S_i\neq {{}0}\} }\\  \leq n\PPP{\xi(0)>u_n, \exists s\in\{-k+1,\dots,k-1\}\setminus\{0\} \mbox{ s.t. }\xi(s)>u_n}\\
 \times \PPP{ \bigcap_{1\leq i\leq k-1}\{S_i\neq 0\}\cap \{S_k=0\}
}.\end{multline*}
According to \eqref{eq:sufficientconditionnotdk}, it follows that

\begin{align}
 \limsup_{n\to\infty}n&\PPP{\left\{\xi({{}0})>u_n, \tilde{M}_{{{}1},{{}k-1}}>u_n, S_{{{}k}}={{}0}\right\} \cap \bigcap_{{{}1}\leq i\leq {{}k-1}}\{S_i\neq {{}0}\} }\nonumber\\
&\qquad\qquad\qquad\qquad\qquad\qquad <\tau\PPP{ 
\bigcap_{1\leq i\leq k-1}\{S_i\neq 0\}\cap \{S_k=0\}
}. \label{eq:maintermdkunpart2}
\end{align} 
Then \eqref{eq:maintermdkun} - \eqref{eq:maintermdkunpart2} implies \eqref{eq:aimdkun}
when $k$ is even. 
 
In a similar way, we prove that the condition $D^{(k)}(u_n)$ is not satisfied when $k$ is {{}odd} by considering this time the event $\{S_{k+1}={{}0}\}$.

\subsection{Proof of Proposition \ref{prop:ppexceedances}}
First, notice that, according to Proposition \ref{prop:existence_extremal_index}, the extremal index of $(\xi(S_n))_{n\geq 0}$ exists and equals $\theta=\sigma q>0$. Now, let $\tau_0$ be such that $p_n^{(\tau_0)}(j)$ converges to some number $p(j)$. In particular $p_n^{(\tau_0)}(j)-p_n^{(\tau_0)}(j-1) \conv[n]{\infty}\theta \pi(j)$ for any $j\geq 1$, where $\pi(j)=\frac{1}{\theta}(p(j)-p(j+1))$. Such a property ensures that Equation (2.5) in \cite{Perfekt94} is satisfied. It follows from Theorem 2.5 in \cite{Perfekt94} that $\Phi_n^{(\tau_0)}$ converges to a point process $\Phi^{(\tau_0)}$ with Laplace transform
\[ L_{\Phi^{(\tau_0)}}(f) := \EEE{\exp\left( -\sum_{x\in \Phi^{(\tau_0)}}f(x)  \right)}  = \exp\left( -\theta \tau_0 \int_0^\infty (1-L(f(t)))\mathrm{d}t  \right),\] for any positive and measurable function $f$, where $L$ denotes the Laplace transform of $\pi$. In particular, $\Phi^{(\tau_0)}$ is a compound Poisson point process of intensity $\theta\tau_0$  with cluster size distribution $\pi$. Moreover, according to Theorem 5.1 in \cite{HHL}, the fact that $\Phi_n^{(\tau_0)}$ converges to $\Phi^{(\tau_0)}$ for some $\tau_0>0$ ensures that  $\Phi_n^{(\tau)}$ converges to $\Phi^{(\tau)}$ for any $\tau>0$. This concludes the proof of Proposition \ref{prop:ppexceedances}.

\section{Examples}
In this section, we give two examples illustrating Proposition \ref{prop:ppexceedances}. The second one extends the first one. Being slightly easier to establish, we have chosen to present Example 1 separately for sake of simplicity.
\label{sec:examples}

\subsection{Example 1}
\label{sec:exiid}
Assume that the $\xi(s)$'s are i.i.d.. Let $N(0)=\#\{i \geq 0: S_i=0\}$ be the number of visits of the random walk $(S_n)_{n\geq 0}$ in site $0$ and recall that \[q=\PPP{S_i\neq 0, \; \forall i\geq 1}=\PPP{N(0)=0}.\]
Notice that $q \in (0,1)$ since the random walk is transient ($p\neq 1/2$) and that $N(0)$ has a geometric distribution with parameter $q$. Moreover, since the $\xi(s)$'s are i.i.d., the extremal index of $(\xi(s))_{s\in \ZZ}$ equals 1. Thus, according to Proposition \ref{prop:existence_extremal_index}, the extremal index of $(\xi(S_n))_{n\geq 0}$ exists and is equal to $\theta=q$.

Now, let $\Phi$ be a compound Poisson point process in $[0,1]$ with intensity $\theta\tau$ and cluster size distribution \begin{align*}
\pi(j) & =\frac{\PPP{N(0)=j}-\PPP{N(0)=j+1}}{\theta}\\
& = q(1-q)^{j-1}.
\end{align*}

\begin{Prop}
\label{prop:compound_explicit}
For any $\tau>0$, the point process of exceedances $\Phi_n$, as defined in \eqref{eq:defppexceedances}, converges in distribution to $\Phi$.
\end{Prop}

\begin{prooft}{Proposition \ref{prop:compound_explicit}}
Since the $\xi(s)$'s are i.i.d., the $\Delta (u_n)$ and the  $D^{(1)}(u_n)$ conditions clearly  hold for $\left(\xi(s)\right)_{s\in\mathbf{Z}}$. Let us justify that $\left(\xi(S_n)\right)_{n\geq 0}$ also satisfies the $\Delta (u_n)$ condition. 
By looking carefully at the proof of \cite[Theorem 2.2]{dHS97}, one can see that the argument, based on a suitable coupling, adapts {\it verbatim} - without assuming that the  $\xi(s)$'s take their values in a finite set - to prove that $\left(\xi(S_n)\right)_{n\geq 0}$ is $\alpha -$mixing. In particular, the $\Delta (u_n)$ condition holds for $\left(\xi(S_n)\right)_{n\geq 0}$.

According to Proposition  \ref{prop:ppexceedances},
 it is sufficient  to compute the cluster size distribution. First, let $(k_n)_{n\geq 0}$ be a family of integers such that $k_n\conv[n]{\infty}\infty$ and $k_n=o(n)$ and let 
\[\Phi_{k_n} = \left\{ \frac{i}{n}: \xi(S_i)>u_n, 0\leq i\leq k_n  \right\}.\]
For any $j\geq 1$, recall that
\[p_n(j) = \PPP{\#\Phi_{k_n} = j | \xi(0)>u_n}.\]
We proceed into two steps: first, we compute the limit of $p_n(j)$; then we compute $\pi(j)$.

\textit{Step 1.} We write
\[p_n(j) = \PPP{\#\Phi_{k_n} \geq  j | \xi(0)>u_n} - \PPP{\#\Phi_{k_n} \geq  j+1 | \xi(0)>u_n}.\]
Moreover, 
\begin{multline*}\PPP{\#\Phi_{k_n} \geq  j | \xi(0)>u_n} = \PPP{\#\Phi_{k_n} \geq j, N_{k_n}(0)\geq j | \xi(0)>u_n} \\+ \PPP{\#\Phi_{k_n} \geq j, N_{k_n}(0)<j | \xi(0)>u_n},
\end{multline*}
where $N_{k_n}(0) = \#\{0\leq i \leq k_n: S_i=0\}$ is the number of visits in 0 until time $k_n$. We show below that the first term of the right-hand side converges to $\PPP{N(0)=j}$ and that the second term converges to 0. First, we notice that, conditional on the event $\{\xi(0)>u_n\}$, we have 
\[\{\#\Phi_{k_n}\geq j\}\cap \{N_{k_n}(0)\geq j\} = \{N_{k_n}(0)\geq j\}.\]
Therefore
\begin{align}
\label{eq:dominant_term_p_n}
\PPP{\#\Phi_{k_n} \geq j, N_{k_n}(0)\geq j | \xi(0)>u_n} & = \PPP{N_{k_n}(0)\geq j|\xi(0)>u_n}\notag\\
& = \PPP{N_{k_n}(0)\geq j}\notag\\
& \conv[n]{\infty}\PPP{N(0)\geq j},
\end{align}
where the second line comes from the fact that the random walk is independent of the scenery and where the third one follows from the Lebesgue's dominated convergence theorem (with $k_n\conv[n]{\infty}\infty$). Moreover, writing $\mathcal{S}_{k_n}=\{S_0,\ldots, S_{k_n}\}$, we have
\begin{align*}
\PPP{\#\Phi_{k_n} \geq j, N_{k_n}(0)<j | \xi(0)>u_n} & \leq \PPP{\exists s \in \mathcal{S}_{k_n}\setminus\{0\}: \xi(s)>u_n|\xi(0)>u_n}\\
& = \frac{\PPP{\{\exists s \in \mathcal{S}_{k_n}\setminus\{0\}: \xi(s)>u_n\}\cap \{\xi(0)>u_n\}}}{\PPP{\xi(0)>u_n}}\\
& \leq \frac{\EEE{\sum_{s\in \mathcal{S}_{k_n}\setminus\{0\}}\PPP{\xi(s)>u_n, \xi(0)>u_n}  }}{\PPP{\xi(0)>u_n}}\\
& \leq \PPP{\xi(0)>u_n}\EEE{\#\mathcal{S}_{k_n}},
\end{align*}
where the last line comes from the fact that the $\xi(s)$'s are i.i.d.. Moreover, according to \cite{LeGRo}, we know that $\frac{\#\mathcal{S}_{k_n}}{k_n}\conv[n]{\infty}q$ a.s.. Thus, according to the Lebesgue's dominated convergence theorem (which can be applied since $\#\mathcal{S}_{k_n}\leq k_n$), we have $\EEE{\frac{\#\mathcal{S}_{k_n}}{k_n}}\conv[n]{\infty}q$. Because $\PPP{\xi(0)>u_n}\eq[n]{\infty} \frac{\tau}{n}$, we have 
\[\PPP{\#\Phi_{k_n} \geq j, N_{k_n}(0)<j | \xi(0)>u_n} = O\left( k_n/n  \right),\]
which converges to 0 since $k_n=o(n)$. This together with \eqref{eq:dominant_term_p_n} gives that 
\[\PPP{\#\Phi_{k_n} \geq j | \xi(0)>u_n} \conv[n]{\infty}\PPP{N(0)\geq j}\] and consequently that 
\[p_n(j) \conv[n]{\infty} \PPP{N(0)\geq j}-\PPP{N(0)\geq j+1}=\PPP{N(0)=j}=:p(j).\]

\textit{Step 2.} We provide below an explicit formula for $\pi$. According to Theorem 4.1 in \cite{Rootzen88} (see also Theorem 2.5 in \cite{Perfekt94}),  we have
\[p(j) = \theta \sum_{m=j}^\infty \pi(m).\]
Therefore
\[\pi(j) = \frac{p(j)-p(j+1)}{\theta} = \frac{\PPP{N(0)=j}-\PPP{N(0)=j+1}}{\theta}.\]
\end{prooft}

\begin{Rk}
It is known that, under suitable conditions, the extremal index can be interpreted as the reciprocal of the mean size of a cluster of exceedances, i.e. $\theta^{-1}=\sum_{j=1}^\infty j\pi(j)$, see e.g. \cite{HHL}. Such an identity holds in the above example since $\frac{1}{q}=\sum_{j=1}^\infty jq(1-q)^{j-1}$.
\end{Rk}

\subsection{Example 2}
Let $k\geq 0$. Assume that $\xi(s)=\max\{ Y_s, Y_{s+1}, \ldots, Y_{s+k} \}$, where the family of random variables $(Y(s))_{s\in \ZZ}$ is assumed to be i.i.d.. For any subset $A\subset\ZZ$, let $N(A)$ be the number of visits of the random walk in $A$, i.e.
\[N(A) = \#\{i\geq 0: S_i\in A\}.\] The following proposition gives an explicit formula for the extremal index of $(\xi(S_n))_{n\geq 0}$ and for the cluster size distribution.

\begin{Prop}
\label{prop:exkdependence}
The point process of exceedances $\Phi_n$, as defined in \eqref{eq:defppexceedances}, converges in distribution to a compound Poisson point process of intensity $\theta \tau$, with $\theta=\frac{q}{k+1}$, and cluster size distribution 
\[\pi(j) = \frac{1}{q} \sum_{A\ni 0: |A|=k+1}\left( \PPP{N(A)=j} - \PPP{N(A)=j+1}  \right),\] for any $\tau>0$. 
\end{Prop}

\begin{prooft}{Proposition \ref{prop:exkdependence}}
Replacing the use of \cite[Theorem 2.2]{dHS97} and its proof by the one of \cite[Theorem 2]{dHKSS03}, one can show that $(\xi(S_n))_{n\geq 0}$ satisfies the $\Delta(u_n)$ condition. Moreover, we notice that the extremal index of $(\xi(s))_{s\in \ZZ}$ is equal to $\sigma=\frac{1}{k+1}$. Moreover, because the $\xi(s)$'s are $k$-dependent, the sequence $(\xi(s))_{s\in \ZZ}$ satisfies the $D^{(k)}(u_n)$. Therefore, according to Proposition \ref{prop:existence_extremal_index}, the extremal index of $(\xi(S_n))_{n\geq 0}$ exists and is equal to $\theta = \frac{q}{k+1}$.

Let us compute the cluster size distribution. We only deal with the case $k=1$ since the general case can be dealt in a similar way. To do it, we  proceed in the same spirit as in Section \ref{sec:exiid}.

\textit{Step 1.} First, we make explicit the limit of
\[p_n(j) = \PPP{\#\Phi_{k_n} \geq  j | \xi(0)>u_n} - \PPP{\#\Phi_{k_n} \geq  j+1 | \xi(0)>u_n}.\]
We notice that
\begin{multline}
\label{eq:clustersizekdependent}
\PPP{\#\Phi_{k_n} \geq  j | \xi(0)>u_n}  = \frac{\PPP{\#\Phi_{k_n} \geq  j, Y_0>u_n}}{\PPP{\xi(0)>u_n}} +  \frac{\PPP{\#\Phi_{k_n} \geq  j, Y_1>u_n}}{\PPP{\xi(0)>u_n}}\\
- \frac{\PPP{\#\Phi_{k_n} \geq  j, Y_0>u_n, Y_1>u_n}}{\PPP{\xi(0)>u_n}}. 
\end{multline}
To deal with the first term of the right-hand side, we write
\[\frac{\PPP{\#\Phi_{k_n} \geq  j, Y_0>u_n}}{\PPP{\xi(0)>u_n}}  = \PPP{\#\Phi_{k_n} \geq  j|Y_0>u_n}\times \frac{\PPP{Y_0>u_n}}{\PPP{\xi(0)>u_n}}.\]
Proceeding in the same spirit as in Section \ref{sec:exiid}, we can prove that $\PPP{\#\Phi_{k_n} \geq  j|Y_0>u_n}$ converges to $\PPP{N(\{-1,0\}) \geq j}$. Moreover, because $\xi(0)=\max\{Y_0,Y_1\}$, where $Y_0$ and $Y_1$ are independent, and because $n\PPP{\xi(0)>u_n}\conv[n]{\infty}\tau$, we have $n\PPP{Y_0>u_n}\conv[n]{\infty}\frac{\tau}{2}$. Therefore, 
\[ \frac{\PPP{\#\Phi_{k_n} \geq  j, Y_0>u_n}}{\PPP{\xi(0)>u_n}} \conv[n]{\infty} \frac{1}{2}\PPP{N(\{-1,0\}) \geq j}. \] In a similar way, we get
 \[ \frac{\PPP{\#\Phi_{k_n} \geq  j, Y_1>u_n}}{\PPP{\xi(0)>u_n}} \conv[n]{\infty} \frac{1}{2}\PPP{N(\{0,1\}) \geq j}. \] 
 Moreover, for the last term of \eqref{eq:clustersizekdependent}, we have
 \[\frac{\PPP{\#\Phi_{k_n} \geq  j, Y_0>u_n, Y_1>u_n}}{\PPP{\xi(0)>u_n}} \leq \frac{(\PPP{Y_0>u_n})^2}{\PPP{\xi(0)>u_n}} \eq[n]{\infty} \frac{\tau}{4n}.\]
 This, together with \eqref{eq:clustersizekdependent}, gives
 \[p_n(j)\conv[n]{\infty} p(j),\]
 with 
 \[p(j) = \frac{1}{2}\left( \PPP{N(\{-1,0\}) = j} +  \PPP{N(\{0,1\}) = j}  \right).\]	

\textit{Step 2.} 
In  the same spirit as in Section \ref{sec:exiid}, we have $\pi(j) = \frac{p(j) - p(j+1)}{\theta}$, that is, since $q=2\theta$
\begin{multline*}
\pi(j) = \frac{1}{q}\Big( \PPP{N(\{-1,0\}) = j}-\PPP{N(\{-1,0\}) = j+1} \\+  \PPP{N(\{0,1\}) = j}-\PPP{N(\{0,1\}) = j+1}\Big).
\end{multline*}
\end{prooft}

\section{Discussion}\label{sec:discussion}

In the examples provided in Section \ref{sec:examples}, the field $\left(\xi(s)\right)_{s\in\mathbf{Z}}$ has a weak dependence property and $\left(\xi(S_n)\right)_{n\geq 0}$ gets stronger mixing properties than the one required in order to apply Propositions \ref{prop:existence_extremal_index} and \ref{prop:ppexceedances}. One may wonder how to proceed to verify that $\left(\xi(S_n)\right)_{n\geq 0}$ satisfies the $\Delta (u_n)$ condition if the field $\left(\xi(s)\right)_{s\in\mathbf{Z}}$ has a stronger dependence property.
We obtain the following partial result in this direction.
\begin{Prop}
\label{prop:mixing}
Let $(S_n)_{n\geq 0}$ be a simple transient random walk, i.e. $p\neq 1/2$ and let $(u_n)_{n\geq 0}$ be a sequence of positive integers. Assume that the following conditions hold.
\begin{enumerate}[(i)]
\item The sequence $(\xi(s))_{s\in \ZZ}$ satisfies the $\Delta(u_n)$ condition. 
\item $\max_{A\in \tilde{\mathcal{B}}_0^n(u_n), B\in \tilde{\mathcal{B}}_0^n(u_n) } \left\{    \left|  \EEE{ \PPP{A | \mathcal{S}_{1:n}} \PPP{B|\mathcal{S}_{1:n}}   } - \PPP{A}\PPP{B} \right| \right\} \conv[n]{\infty} 0.$
\end{enumerate}
Then $(\xi(S_n))_{n\geq 0}$ satisfies the $\Delta(u_n)$ condition. 
\end{Prop}

\begin{Rk}
In the above proposition, the notation $\tilde{\mathcal{B}}_0^n(u_n)$ stands for the $\sigma$-algebra generated by events of the form $\{\xi(S_i) \leq u_n\}$, $0\leq i\leq n$, and   $\mathcal{S}_{0:n}=\{S_0, \ldots, S_n\}$. The assumption \textit{(ii)} is probably unnecessary, our result should hold if we only assume that $(\xi(s))_{s\in \ZZ} $ satisfies the $\Delta(u_n)$ condition. At least assumption \textit{(ii)} can be replaced by the slightly weaker assumption \textit{(ii')} stated below. Assumption \textit{(ii')} seems more natural from the point of view of EVT but is a bit heavier to state. Moreover, it appears in our examples that it is not easier to check the latter than assumption \textit{(ii)}.

\begin{enumerate}[\textit{(ii')}]\item For some $\beta \in (1/2,1)$, 
\[  \max_{1\leq k\leq n-\tilde{\ell}_n} \max_{A\in \tilde{\mathcal{B}}_0^k(u_n), B\in \tilde{\mathcal{B}}_{k+\tilde{\ell}_n}^n(u_n) } \left\{    \left|  \EEE{ \PPP{A | \mathcal{S}_{0:k}} \PPP{B|\mathcal{S}_{k+\tilde{\ell}_n:n}}   } - \PPP{A}\PPP{B} \right| \right\} \conv[n]{\infty} 0,\]
where \begin{equation}
\label{eq:defln}
 \tilde{\ell}_n = \left\lceil  \frac{\ell_{2n+1} + 2n^\beta}{|2p-1|} \right\rceil.
 \end{equation} 
 \end{enumerate} 
\end{Rk}

\begin{prooft}{Proposition \ref{prop:mixing}}
Without loss of generality, we only deal with the case $p>1/2$.  Let $(\ell_n)_{n\geq 1}$ be such that  $\alpha_{n,\ell_n}\conv[n]{\infty}0$ and  $\ell_n=o(n)$. Let $\tilde{\ell}_n$ be as in \eqref{eq:defln}, 
with $\beta \in (1/2,1)$. Notice that $\tilde{\ell}_n\conv[n]{\infty}\infty$ and $\tilde{\ell}_n = o(n)$. Similarly to $\alpha_{n,\ell}$, we introduce the term
 \[\tilde{\alpha}_{n,\ell} =  \max_{1\leq k\leq n-\ell} \max_{A\in \tilde{\mathcal{B}}_0^k(u_n), B\in \tilde{\mathcal{B}}_{k+\ell}^n(u_n) }\{ |\PPP{A\cap B} - \PPP{A}\PPP{B}|\},\] where
  $\tilde{\mathcal{B}}_{m_1}^{m_2}(u_n)$ denotes the $\sigma$-field generated by events $\{\xi(S_i)\leq u_n\}$, $m_1\leq i\leq m_2$. We prove below that $\tilde{\alpha}_{\tilde{\ell}_n,n}\conv[n]{\infty}0$. To do it, we will use the following lemma whose proof is given in Appendix \ref{app:proofconcentrationinequality}.
  
  \begin{Le}
  \label{Le:concentrationinequality}
Assume that $p>1/2$. Let us consider the following event:
  \[E_n= \bigcap_{k\leq n-\tilde{\ell}_n}\left\{ \max\{S_0,\ldots, S_k\} \leq k(2p-1)+n^\beta    \right\}\cap \left\{  \min\{S_{k+\tilde{\ell}_n}, \ldots, S_n\} \geq (k+\tilde{\ell}_n)(2p-1) - n^\beta \right\}    .\]
  Then there exist two positive constants $c_1$ and $c_2$ such that $\PPP{E_n^c} \leq c_1 e^{-c_2n^{2\beta-1}}$.
  \end{Le}
In particular, Lemma \ref{Le:concentrationinequality} ensures that the event $E_n$ occurs with high probability since $\beta>1/2$. 

Now, let $k\leq n-\tilde{\ell}_n$ be fixed.  Let us consider two events $A\in \tilde{\mathcal{B}}_0^k(u_n)$ and $B\in \tilde{\mathcal{B}}_{k+\tilde{\ell}_n}^n(u_n)$. First, according to Lemma \ref{Le:concentrationinequality}, we notice that
\begin{align*}
\PPP{A\cap B} & = \PPP{A\cap B\cap E_n} + o(1)\\
& = \EEE{\ind{E_n} \PPP{A\cap B | \mathcal{S}_{1:n}}  } + o(1),
\end{align*}
where $o(1)$ only depends on $n$. Moreover, conditional on $\mathcal{S}_{1:n}$ and on the event $E_n$, we know that 
\begin{equation} \label{eq:conditionalindependence} \left| \ind{E_n} \PPP{A\cap B | \mathcal{S}_{1:n}}  - \ind{E_n} \PPP{A| \mathcal{S}_{1:n}} \PPP{B| \mathcal{S}_{1:n}}   \right| \leq \alpha_{2n+1, \ell_{2n+1}}.  \end{equation}
Indeed, since $A\in \tilde{\mathcal{B}}_0^k(u_n)$ and since we are on $E_n$, conditional on $\mathcal{S}_{1:n}$, the event $A$ only depends of events of the form $\{\xi(s) \leq u_n\}$, with $s\leq k(2p-1)+n^\beta$. In the same way, the event $B$ only depends of events of the form $\{\xi(s) \leq u_n\}$, with $s\geq (k+\tilde{\ell}_n)(2p-1)-n^\beta$. Equation  \eqref{eq:conditionalindependence} follows since the difference between $k(2p-1)+n^\beta$ and $(k+\tilde{\ell}_n)(2p-1)-n^\beta$ is at least $\ell_{2n+1}$. Now, because $(\xi(s))_{s\in \ZZ}$ satisfies the $\Delta(u_n)$ condition, we know that   $\alpha_{2n+1, \ell_{2n+1}}$ converges to 0 as $n$ goes to infinity. Thus, thanks again to Lemma \ref{Le:concentrationinequality}, we have
\[ \PPP{A\cap B}  = \EEE{\PPP{A| \mathcal{S}_{1:n}} \PPP{B| \mathcal{S}_{1:n}}  }  + o(1).\]
This together with assumption \textit{(ii)} concludes the proof of Proposition \ref{prop:mixing}. 
\end{prooft}

\appendix
\section{Proof of Lemma \ref{Le:concentrationinequality}}\label{app:proofconcentrationinequality}

It is sufficient to prove that the events 
\[\bigcup_{k\leq n-\tilde{\ell}_n}\left\{ \max\{S_0,\ldots, S_k\} > k(2p-1)+n^\beta    \right\}\]
and 
\[\bigcup_{k\leq n-\tilde{\ell}_n} \left\{  \min\{S_{k+\tilde{\ell}_n}, \ldots, S_n\} < (k+\tilde{\ell}_n)(2p-1) - n^\beta \right\}    \]
occur with probability smaller than $c_1 e^{-c_2n^{2\beta-1}}$. We only deal with the first one since the second one can be dealt in a similar way. To do it, we notice that
\begin{align*}
\PPP{\bigcup_{k\leq n-\tilde{\ell}_n}\left\{ \max\{S_0,\ldots, S_k\} > k(2p-1)+n^\beta    \right\}} & \leq \sum_{k\leq n-\tilde{\ell}_n}\sum_{i\leq k}\PPP{S_i>k(2p-1)+n^\beta }\\
& \leq \sum_{k\leq n-\tilde{\ell}_n} k \PPP{S_k>k(2p-1)+n^\beta },
\end{align*}
where the last line comes from the fact that $S_i$ is stochastically dominated by $S_k$, with $i\leq k$, since $p>1/2$. Now, let $k\leq  n-\tilde{\ell}_n$ be fixed, and let $S'_k$ be a binomial random variable with parameter $(k,p)$, so that $S_k \overset{sto}{=}2S'_k-k$. We have
\[\PPP{S_k > k(2p-1) + n^\beta}  = \PPP{S'_k>kp+\frac{1}{2}n^\beta}  \leq  e^{-\tfrac{1}{2} \cdot\tfrac{n^{2\beta}}{k}},\]
according to the Hoeffding's inequality. Since $k\leq n$, the last term is lower than $e^{-\tfrac{1}{2}n^{2\beta-1}}$. Summing over $k$, we get
\[\PPP{\bigcup_{k\leq n-\tilde{\ell}_n}\left\{ \max\{S_0,\ldots, S_k\} > k(2p-1)+n^\beta    \right\}} \leq n^2 e^{-\tfrac{1}{2}n^{2\beta-1}}.\] This concludes the proof of Lemma \ref{Le:concentrationinequality}.


\begin{thebibliography}{99}

\bibitem{Chen_Dar}Chenavier, N. \& Darwiche, A. Extremes for transient random walks in random sceneries under weak independence conditions. {\em Statist. Probab. Lett.}. \textbf{158} pp. 108657, 6 (2020).

\bibitem{Dar2}Chenavier, N. \& Darwiche, A. Some properties on extremes for transient random walks in random sceneries. {\em Available In Https://arxiv.org/pdf/2210.04854.pdf} (2022)

\bibitem{CHM}Chernick, M., Hsing, T. \& McCormick, W. Calculating the extremal index for a class of stationary sequences. {\em Adv. In Appl. Probab.}. \textbf{23}, 835-850 (1991).

\bibitem{dHKSS03}Hollander, F., Keane, M., Serafin, J. \& Steif, J. Weak Bernoullicity of random walk in random scenery. {\em Japan. J. Math. (N.S.)}. \textbf{29}, 389-406 (2003).

\bibitem{dHS97}Hollander, F. \& Steif, J. Mixing properties of the generalized $T,T^{-1}$-process. {\em J. Anal. Math.}. \textbf{72} pp. 165-202 (1997).

\bibitem{EKM}Embrechts, P., Klüppelberg, C. \& Mikosch, T. Modelling extremal events, volume 33 of Applications of Mathematics (New York). {\em Springer-Verlag}, Berlin,  For insurance and finance (1997).

\bibitem{franke_saigo2009bis}Franke, B. \& Saigo, T. The extremes of a random scenery as seen by a random walk in a random environment. {\em Statist. Probab. Lett.}. \textbf{79}, 1025-1030 (2009).

\bibitem{franke_saigo_2009}Franke, B. \& Saigo, T. The extremes of random walks in random sceneries. {\em Advances In Applied Probability}. \textbf{41}, 452-468 (2009)

\bibitem{HHL}Hsing, T., Hüsler, J. \& Leadbetter, M. On the exceedance point process for a stationary sequence. {\em Probab. Theory Related Fields}. \textbf{78}, 97-112 (1988).

\bibitem{KPN}Katz, R., Parlange, M. \& Naveau, P. Statistics of extremes in hydrology. {\em Advances In Water Resources}. \textbf{25} pp. 1287-1304 (2002).

\bibitem{KS}Kesten, H. \& Spitzer, F. A limit theorem related to a new class of sel-similar processes. {\em Z. Wahrsch. Verw. Gebiete}. \textbf{50} pp. 5-25 (1979).

\bibitem{LeGRo}Le Gall, J. \& Rosen, J. The range of stable random walks. {\em Ann Proba}. \textbf{19}, 650-705 (1991). 

\bibitem{L2}Leadbetter, M. Extremes and local dependence in stationary sequences. {\em Z. Wahrsch. Verw. Gebiete}. \textbf{65}, 291-306 (1983).

\bibitem{LLR}Leadbetter, M., Lindgren, G. \& Rootzén, H. Extremes and related properties of random sequences and processes. {\em Springer-Verlag} (1983).

\bibitem{Perfekt94}Perfekt, R. Extremal behaviour of stationary Markov chains with applications. {\em Ann. Appl. Probab.}. \textbf{4}, 529-548 (1994).

\bibitem{Rootzen88}Rootzén, H. Maxima and exceedances of stationary Markov chains. {\em Adv. In Appl. Probab.}. \textbf{20}, 371-390 (1988).

\bibitem{YGLN}Yiou, P., Goubanova, K., Li, Z. \& Nogaj, M.  Weather regime dependence of extreme value statistics for summer temperature and precipitation. {\em Nonlinear Processes In Geophysics}. \textbf{15}, 365-378 (2008).




\end{thebibliography}

\end{document}